\title{\vspace{-0.5cm}How many random edges make a dense hypergraph
non-$2$-colorable?}

\author{
Benny Sudakov \thanks{Department of Mathematics,
Princeton University, Princeton, NJ 08544. E-mail:
bsudakov@math.princeton.edu.
Research supported in part by NSF
grant DMS-0355497, USA-Israeli BSF grant, and by an Alfred P. Sloan
fellowship.}
\and Jan Vondr\'ak \thanks{ Microsoft Research, Redmond, WA
98502. E-mail: vondrak@microsoft.com.}
}

\documentclass [11pt]{article}
\usepackage{amsfonts}
\usepackage{amsmath}
\usepackage{latexsym}
\usepackage{epsfig}
\usepackage{graphicx}

\newtheorem{theorem}{Theorem}[section]
\newtheorem{lemma}[theorem]{Lemma}
\newtheorem{proposition}[theorem]{Proposition}

\newtheorem{remark}[theorem]{Remark}

\def\HL{H}     % notation for the lower-bound example

\date{}
\begin{document}
\maketitle
\begin{abstract}
We study a model of random uniform hypergraphs, where a random
instance is obtained by adding random edges to a large hypergraph
of a given density. The research on this model for graphs has been
started by Bohman et al. in \cite{BFM}, and continued in
\cite{BFKM} and \cite{KST}. Here we obtain a tight bound on the
number of random edges required to ensure non-$2$-colorability. We
prove that for any $k$-uniform hypergraph with
$\Omega(n^{k-\epsilon})$ edges, adding $\omega(n^{k\epsilon/2})$
random edges makes the hypergraph almost surely non-$2$-colorable.
This is essentially tight, since there is a $2$-colorable
hypergraph with $\Omega(n^{k-\epsilon})$ edges which almost surely
remains $2$-colorable even after adding $o(n^{k \epsilon / 2})$
random edges.
\end{abstract}

\section{Introduction}

Research on random graphs and hypergraphs has a long history with
thousands of papers and two monographs by Bollob\'as \cite{B} and
by Janson et al. \cite{JLR} devoted to the subject and its diverse
applications. In the classical Erd\H{o}s-R\'{e}nyi model
\cite{ER}, a random graph is generated by starting from an empty
graph and then adding certain number of random edges. More
recently, Bohman, Frieze and Martin \cite{BFM} considered a
generalized model where one starts with a fixed graph $G=(V,E)$
and then inserts a collection $R$ of additional random edges. We
denote the resulting random graph by $G+R$. The initial graph $G$
can be regarded as given by an adversary, while the random
perturbation $R$ represents noise or uncertainty, independent of
the initial choice. This scenario is analogous to the {\em
smoothed analysis} of algorithms proposed by Spielman and Teng
\cite{ST}, where an algorithm is assumed to run on the worst-case
input, modified by a small random perturbation.

Usually, one investigates {\em monotone properties} of random
graphs or hypergraphs; i.e., properties which cannot be destroyed by
adding more edges, like the property of containing a certain fixed
subgraph. Given a monotone property $\cal A$ of graphs on $n$ vertices,
we can ask what are the parameters for which a random
graph has property $\cal A$ almost surely, i.e. with probability
tending to $1$ as the number of vertices $n$ tends to infinity.
In our setting, we start with a  fixed hypergraph $H$ and inquire
how many random edges $R$ we have to add so that $H+R$ has
property $\cal A$ almost surely. This question is too general
to get concrete and meaningful results, valid for all hypergraphs $H$.
Therefore, rather than considering a completely arbitrary
$H$, we start with a hypergraph from a certain natural class.
One such class of graphs was considered in
\cite{BFM}, where the authors analyze the question of how many
random edges need to be added to a graph $G$ of minimal degree at least $dn, 0<d<1$,
so that the resulting graph $G+R$ is almost
surely Hamiltonian. Further properties of random graphs in this model are
explored in \cite{BFKM}.

In \cite{KST}, Krivelevich et al. considered a slightly more
general setting, in which one performs a small random perturbation
of a graph $G$ with at least $dn^2$ edges. Observe that since $G$
has at least $dn^2$ edges, removing a small set of random edges
would leave the total number of edges in $G$ essentially
unchanged. Therefore one only has to focus on the case of adding
random edges. In \cite{KST}, the authors obtained tight results
for the appearance of a fixed subgraph and for certain Ramsey
properties in this model. In the same paper, they also considered
random formulas obtained by adding random $k$-clauses
(disjunctions of $k$ literals) to a fixed $k$-SAT formula.
Krivelevich et al. proved that for any formula with at least $n^{k
- \epsilon}$ $k$-clauses, adding $\omega(n^{k \epsilon})$ random
clauses of size $k$ makes the formula almost surely unsatisfiable.
This is tight, since there is a $k$-SAT formula with
$n^{k-\epsilon}$ clauses which almost surely remains satisfiable
after adding $o(n^{k\epsilon})$ random clauses. A related
question, which was raised in \cite{KST}, is to find a threshold
for non-$2$-colorability of a random hypergraph obtained by adding
random edges to a large hypergraph of a given density.

For an integer $k \geq 2$, a {\em $k$-uniform hypergraph}
is an ordered pair $H=(V,E)$, where $V$ is a finite non-empty set, called set of
{\em vertices} and $E$ is a family of distinct $k$-subsets of $V$, called
the {\em edges} of $H$. A {\em $2$-coloring} of a hypergraph $H$ is a partition of
its vertex set $V$ into two color classes so that no edge in $E$ is monochromatic.
A hypergraph which admits a $2$-coloring is called {\em $2$-colorable}.

$2$-colorability is one of the fundamental properties of
hypergraphs, which was first introduced and studied by Bernstein
\cite{Ber} in 1908 for infinite hypergraphs. $2$-colorability in
the finite setting, also known as ``Property B" (a term coined by
Erd\H{o}s in reference to Bernstein), has been studied extensively
in the last forty years (see, e.g., \cite{Er, Er1, EL, Beck, RS}).
While $2$-colorability of graphs is well understood being
equivalent to non-existence of odd cycles, for $k$-uniform
hypergraphs with $k \geq 3$ it is already $NP$-complete to decide
whether a $2$-coloring exists \cite{L}. Consequently, there is no
efficient characterization of $2$-colorable hypergraphs. The
problem of $2$-colorability of random $k$-uniform hypergraphs for
$k \geq 3$ was first studied by Alon and Spencer \cite{AS}. They
proved that such hypergraphs with $m=(c2^k/k^2)n$ edges are almost
surely $2$-colorable. This bound was improved later by Achlioptas
et al. \cite{AKKT}. Recently, the threshold for $2$-colorability
has been determined very precisely. In \cite{AM} it was proved
that the number of edges for which a random $k$-uniform hypergraph
becomes almost surely non-$2$-colorable is $(2^{k-1} \ln 2 - O(1))
n$.

Interestingly, the threshold for non-$2$-colorability is roughly
one half of the threshold for $k$-SAT. It has been shown in
\cite{AP} that a formula with $m$ random $k$-clauses becomes
almost surely unsatisfiable for $m = (2^k \ln 2 - O(k)) n$. The
two problems seem to be intimately related and it is natural to
ask what is their relationship in the case of a random
perturbation of a fixed instance. Recall that from \cite{KST} we
know that for any $k$-SAT formula with $n^{k-\epsilon}$ clauses,
adding $\omega(n^{k \epsilon})$ random clauses makes it almost
surely unsatisfiable. In fact, the same proof yields that for any
$k$-uniform hypergraph $H$ with $n^{k-\epsilon}$ edges, adding
$\omega(n^{k \epsilon})$ random edges destroys $2$-colorability
almost surely. Nonetheless, it turns out that this is not the
right answer. It is enough to use substantially fewer random edges
to destroy $2$-colorability: roughly a square root of the number
of random clauses necessary to destroy satisfiability. The
following is our main result.

\begin{theorem}
\label{main} Let $k,\ell \geq 2$, $\epsilon \geq 0$ be fixed and
let $H$ be a $2$-colorable $k$-uniform hypergraph with
$\Omega(n^{k-\epsilon})$ edges. Then the hypergraph $H'$ obtained
by adding to $H$ a collection $R$ of $\omega\big(n^{\ell \epsilon
/ 2}\big)$ random $\ell$-tuples is almost surely
non-$2$-colorable.
\end{theorem}

Observe that for $\epsilon \geq 2 / \ell$, the result is easy.
Regardless of the hypergraph $H$, it is well known that a
collection of $\omega(n)$ random $\ell$-tuples on $n$ vertices is
almost surely non-$2$-colorable. So we will be only interested in
the case when $\epsilon < 2 / \ell$. For such $\epsilon$ we obtain
the following result, which shows that the assertion of Theorem
\ref{main} is essentially best possible.

\begin{theorem}
\label{main-construction} For fixed $k,\ell \geq 2$ and $0 \leq
\epsilon < 2 / \ell$, there exists a $2$-colorable $k$-uniform
hypergraph $H$ with $\Omega(n^{k-\epsilon})$ edges such that its
union with a collection $R$ of $o\big(n^{\ell \epsilon / 2}\big)$
random $\ell$-tuples is almost surely $2$-colorable.
\end{theorem}

The rest of this paper is organized as follows. In the next section we
present an example of the hypergraph which shows that our main result is
essentially best possible. In Section 3 we discuss some natural difficulties in
proving Theorem \ref{main} and describe how to deal with them in the case of
bipartite graphs. This result also serves as a basis for induction which we use
in Section 4 to prove the general case of $2$-colorable $k$-uniform hypergraphs.

\begin{remark}
\label{equiv-model} We have two alternative ways of adding random
edges. Either we can sample a random $\ell$-tuple $|R|$ times,
each time uniformly and independently from the set of all ${n
\choose \ell}$ $\ell$-tuples. Or we can pick each $\ell$-tuple
randomly and independently with probability $p=|R|/{n \choose
\ell}$. Since $2$-colorability is a monotone property, it follows,
as in Bollob\'as \cite{B}, Theorem 2.2 and a similar remark in
\cite{KST}, that if the resulting hypergraph is almost surely
non-$2$-colorable ($2$-colorable) in one model then this is true
in the other model as well. This observation can sometimes
simplify our calculations.
\end{remark}

\paragraph{\bf Notation.} Let $H = (V,E)$ be a $k$-uniform hypergraph.
In the following, we use the notions of {\em degree} and {\em
neighborhood}, generalizing their usual meaning in graph theory.
For a vertex $v \in V$, we define its degree $d(v)$ to be the
number of edges of $H$ that contain $v$. More generally, for a
subset of vertices $A \subset V, |A| < k$, we define its degree $
d(A) = |\{e \in E: A \subset e\}|.$ For a $(k-1)$-tuple of
vertices $A$, we define its {\em neighborhood} as $ N(A) = \{w \in
V \setminus A: A \cup \{w\} \in E\}.$  Also, for a $(k-2)$-tuple
of vertices $A$, we define its {\em link} as $ \Gamma(A) =
\{\{u,v\} \in V \setminus A: A \cup \{u,v\} \in E \}.$

Throughout the paper we will systematically omit floor and ceiling
signs for the sake of clarity of presentation. Also, we use the
notations $a_n =\Theta(b_n)$, $a_n =O(b_n)$ or $a_n =\Omega(b_n)$
for $a_n, b_n >0$ and $n \to \infty$ if there are absolute
constants $C_1$ and $C_2$ such that $C_1\, b_n < a_n < C_2\, b_n$,
$a_n < C_2 b_n$ or $a_n > C_1 b_n$ respectively. The notation
$a_n=o(b_n)$ means that $a_n/b_n \rightarrow 0$ as $n \rightarrow
\infty$, and $a_n= \omega(b_n)$ means $a_n/b_n \rightarrow
\infty$. The parameters $k,\ell,\epsilon$ are considered constant.

\section{The lower bound}
The following example proves Theorem \ref{main-construction} and
shows that our main result is essentially best possible.

\vspace{-0.2cm}
\paragraph{\bf Construction.} Partition the set of vertices $[n]$ into
three disjoint subsets $X, Y, Z$ where $|X|=|Y|=n^{1-\epsilon/2}$.
Let $\HL$ be a $k$-uniform hypergraph whose edge set consists of
all $k$-tuples which have exactly one vertex in $X$, one vertex in
$Y$ and $k-2$ vertices in $Z$. By definition the number of edges
in $\HL$ is $\Theta(n^{k-\epsilon})$.

\begin{figure}[!ht]
\begin{center}
\scalebox{0.6}{\includegraphics{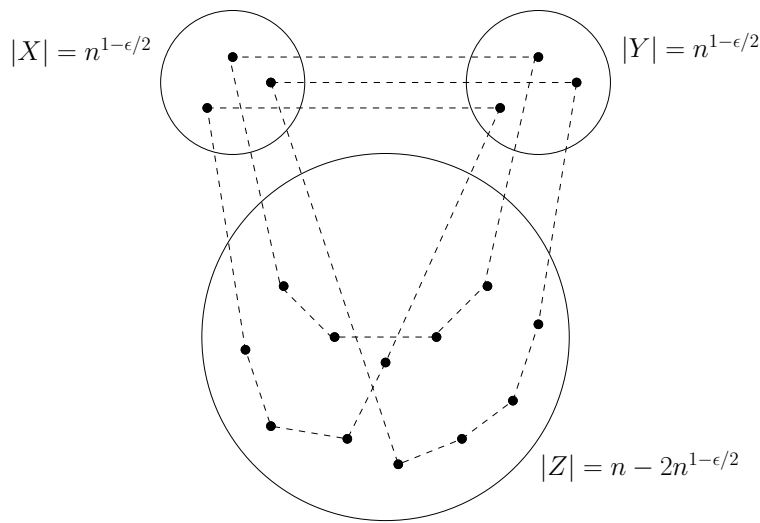}}
\caption{Construction of the hypergraph $\HL$.}
\end{center}
\end{figure}

\vspace{-0.2cm}
\paragraph{\bf Claim.} Color all the vertices in $X$ by
color $1$ and vertices in $Y$ by color $2$. Note that no matter
how we assign colors to the remaining vertices, this gives a
proper $2$-coloring of $\HL$. Let $R$ be a set of $o\big(n^{\ell
\epsilon/2}\big)$ random $\ell$-tuples. Then almost
surely we can $2$-color $Z$ so that none of the $\ell$-tuples in
$R$ is monochromatic, i.e., there exists a proper $2$-coloring of
$\HL+R$.

To prove this claim we transform $R$ into another random instance
$R'$ that contains only single vertices with a fixed {\em
prescribed color} and edges of size two which must not be
monochromatic. Following Remark~\ref{equiv-model} we can assume
that $R$ was obtained by choosing every $\ell$-tuple in $[n]$
randomly and independently with probability
$p=o\left(n^{\ell\epsilon/2-\ell}\right)$. First note that almost
surely there is no $\ell$-tuple in $R$ whose vertices are all in
$X$ or in $Y$. Indeed, since $|X|=|Y|=n^{1-\epsilon/2}$, the
probability that there is such an $\ell$-tuple is at most $2
{n^{1-\epsilon/2} \choose \ell}\,p =o(1)$. Also, every
$\ell$-tuple in $R$ which has vertices in both $X$ and $Y$ is
already $2$-colored so we discard it.

For every $v \in Z$ we add it to $R'$ with prescribed color $1$ if
there is a subset $A$ of $Y$ of size $\ell-1$ such that $A \cup
\{v\} \in R$. Since $\epsilon< 2/\ell \leq 1$, the probability of
this event is
$$p_1={|Y| \choose \ell-1}\,p = {n^{1-\epsilon/2} \choose \ell-1}\, p \leq
n^{(\ell-1)(1-\epsilon/2)}
\,p=o\big(n^{-1+\epsilon/2}\big)=o\big(n^{-1/2}\big).$$ Similarly,
if there is a subset $B$ of $X$ of size $\ell-1$ such that $B \cup
\{v\} \in R$ then we add $v$ to $R'$ with prescribed color $2$.
The probability $p_2$ of this event is also $o\big(n^{-1/2}\big)$.

Fix an ordering $v_1< v_2< \dots$ of all vertices in $Z$. For
every pair of vertices $u, w \in Z$ we add an edge $\{u,w\}$ to
$R'$ if there is an $\ell$-tuple $L \in R$ such that the two
smallest vertices in $L\cap Z$ are $u$ and $w$. Since the number of such
possible $\ell$-tuples is at most ${n \choose \ell-2}$, and
$\epsilon < 2/\ell$, the probability of this event is
$$p_3 \leq {n \choose \ell-2} p =  O\left(n^{\ell-2}p\right) =
 o\left(n^{\ell\epsilon/2-2}\right) = o\left(n^{-1}\right).$$
 Also note that by definition all the above events are independent since they depend
on disjoint sets of $\ell$-tuples. By our construction, any
$2$-coloring of $Z$ in which singletons in $R'$ get prescribed
colors and no $2$-edge is monochromatic gives a proper
$2$-coloring of $R$. Therefore, to complete the proof of Theorem
\ref{main-construction}, it is enough to prove the following
simple statement.

\begin{lemma}
\label{last} Let $R'$ be a random instance which is obtained as
follows. For $i=1,2$ we choose every vertex in $[n]$ with
probability $p_i= o\big(n^{-1/2}\big)$ (independently for $i=1,2$)
and prescribe to it color $i$. In addition we choose every pair of
vertices to be an edge in $R'$ with probability $p_3=o(n^{-1})$.
Then almost surely there exists a $2$-coloring of $[n]$ in which
all singletons in $R'$ get prescribed colors and no edge is
monochromatic.
\end{lemma}

\noindent {\bf Proof.}\, Let $G$ be the graph formed by edges from
$R'$. The probability that there is a vertex with conflicting
prescribed colors is $n p_1 p_2 = o(1)$. The probability that $G$
contains a cycle is at most $\sum_{s=3}^n n^s p_3^s
 = O(n^3 p_3^3) = o(1)$. Finally the probability that there exists a path
between two vertices with any prescribed color is also bounded by
$$ \sum_{s=1}^n {n \choose 2} (p_1+p_2)^2 n^{s-1} p_3^s
= o\big(n (p_1+p_2)^2\big) = o(1).$$

Therefore almost surely no vertex gets prescribed conflicting
colors, every connected component of $G$ is a tree and contains at
most one vertex with prescribed color. This immediately implies
the assertion of the lemma, since every tree can be $2$-colored,
starting from the vertex with prescribed color (if any). \hfill
$\Box$

\section{Bipartite graphs}

Now let's turn to Theorem~\ref{main}. First, consider the case of
$k=\ell=2$. Here, we claim that for any bipartite graph $G$ with
$\Omega\big(n^{2-\epsilon}\big)$ edges, adding
$\omega(n^{\epsilon})$ random edges makes the graph almost surely
non-bipartite. This will follow quite easily, since it turns out
that almost surely we will insert an edge inside one part of a
bipartite connected component of $G$, creating an odd cycle (see
the proof of Proposition \ref{bipartite}).

However, with the more general hypergraph case in mind, we are
also interested in a scenario where random $\ell$-tuples are added
to a bipartite graph, and $\ell > 2$. Then we ask what is the
probability that the resulting hypergraph is $2$-colorable (i.e.,
no $2$-edge and no $\ell$-edge should be monochromatic). We prove
the following special case of Theorem~\ref{main}.

\begin{proposition}
\label{bipartite} Let $\ell \geq 2$, $0 \leq \epsilon < 2/\ell$
and let $G$ be a bipartite graph with
$\Omega\big(n^{2-\epsilon}\big)$ edges. Then the hypergraph $H$
obtained by adding to $G$ a collection $R$ of $\omega\big(n^{\ell
\epsilon/2}\big)$ random $\ell$-tuples is almost surely
non-$2$-colorable.
\end{proposition}

\begin{figure}[!ht]
\label{fig-bipartite}
\begin{center}
\scalebox{0.8}{\includegraphics{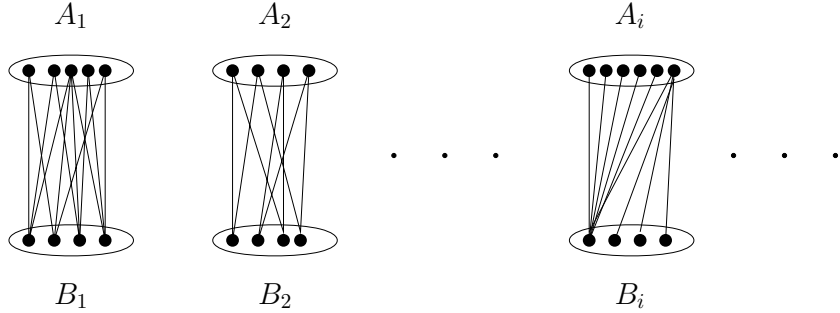}}
\caption{Components of the bipartite graph $G$.}
\end{center}
\end{figure}

\noindent {\bf Proof.}\, Consider the connected components of $G$
which are bipartite graphs on disjoint vertex sets $(A_1, B_1)$,
$(A_2, B_2), \ldots$ (see Figure 2). Denote $a_i = |A_i|$, $b_i =
|B_i|$ and assume $a_i \geq b_i$. The number of edges in each
component is at most $a_i b_i$. Since the total number of edges is
at least $c n^{2 - \epsilon}$ for some constant $c>0$, we have
$$ \sum{a_i^2} \geq \sum{a_i b_i} \geq c n^{2-\epsilon}.$$
Observe that for $\ell=2$, the number of pairs of vertices inside
the sets $\{A_i\}$ is $\sum{a_i \choose 2} \geq \frac12 (c
n^{2-\epsilon} - n) \geq c' n^{2-\epsilon}$, so a random edge
lands inside one of these sets with probability at least $c'
n^{-\epsilon}$. Consequently, the probability that none of the
$\omega(n^\epsilon)$ random edges ends up inside some $A_i$ is at
most $(1 - c' n^{-\epsilon}) ^ {\omega(n^{\epsilon})} = o(1)$.
Thus almost surely, $G + R$ contains an odd cycle.

On the other hand, in the general case we are adding $\omega\big(n^{\ell
\epsilon/2}\big)$ random $\ell$-tuples, which might never end up
inside any vertex set $A_i$. The probability of hitting a specific
$A_i$ is ${a_i \choose \ell} / {n \choose \ell} = O\big(a_i ^ \ell /
n^\ell\big)$. For example, if $G$ has $n^\epsilon$ components with
$a_i = b_i = n^{1-\epsilon}$, then this probability is at most
$O\big(\sum{a_i^\ell} / n^\ell\big) = O\big(n^{-(\ell-1) \epsilon}\big)$. Hence we
need $\omega\big(n^{(\ell-1) \epsilon}\big)$ random $\ell$-tuples to hit
almost surely some $A_i$. This suggests a difficulty with the
attempt to place a random $\ell$-tuple in a set which is forced to
be monochromatic by the original graph. We have to allow ourselves
more freedom and consider sets which are monochromatic only under
certain colorings.

More specifically, each of the sets $A_i$, $B_i$ must be
monochromatic under any coloring, and at least half of them must
share the same color. We do not know a priori which sets will
share the same color, yet we can estimate the probability that
{\em any} of these configurations allows a feasible coloring
together with the random $\ell$-tuples. First, it is convenient to
assume that the sets have roughly equal size, in which case we
have the following claim.

\begin{lemma}
\label{clusters} Suppose we have $t$ disjoint subsets $A_1,
\ldots, A_t$ of $[n]$ of size $\Theta(n^{1-\alpha})$. Let $\alpha
\geq \epsilon/2$, $t = \Omega\big(n^{\frac{\ell}{\ell-1} (\alpha -
\epsilon/2)}\big)$ and let $R$ be a collection of $\omega\big(n^{\ell
\epsilon / 2}\big)$ random $\ell$-tuples on $[n]$. Then the
probability that $R$ can be $2$-colored in such a way that each
$A_i$ is monochromatic is at most $e^{-\omega(t)}$.
\end{lemma}

\noindent {\bf Proof.} Consider the $2^t$ possible colorings in
which all $A_i$ are monochromatic. For each such coloring there is
a set of indices $I, |I|\geq t/2$ such that the sets $A_i, i \in
I$ share the same color. Since $A_i$ are disjoint we have
$|\cup_{i \in I}{A_i}| \geq c_1 t n^{1-\alpha}$ for some $c_1>0$.
The probability that one random $\ell$-tuple falls inside this set
is at least ${c_1 t n^{1-\alpha} \choose \ell} / {n \choose \ell}
\geq c_2 (t n^{-\alpha})^\ell$ for some $c_2 > 0$. Since
$t^{\ell-1} = \Omega\big(n^{\ell (\alpha - \epsilon/2)}\big)$, it
implies that
$$ \Pr\big[\mbox{ $\cup_{i\in I} A_i$ contains no $\ell$-tuple from
$R$}\big]
 \leq \Big(1 - (c_2 t n^{-\alpha})^\ell \Big)^{\omega(n^{\ell \epsilon /
 2})}  \leq e^{-\omega(t^\ell n^{-\ell (\alpha - \epsilon/2)})} =
 e^{-\omega(t)}. $$
Therefore, by the union bound over all choices of $I$, we get
$$ \Pr\big[\exists\, I\, \mbox{such that $\cup_{i\in I} A_i$ contains no $\ell$-tuple from
$R$}\big] \leq 2^t e^{-\omega(t)} = e^{-\omega(t)}.$$
In particular, almost surely there is no $2$-coloring of $R$ in
which all $A_i$ are monochromatic.
 \hfill $\Box$

Now we can finish the proof of Proposition~\ref{bipartite} for
$\ell \geq 3$. Partition the components of $G$ according to their
size and let $G_s$ contain all the components with $|A_i| \in
[2^{s-1},2^s)$. If there is any $A_i$ of size at least $n^{1 -
\epsilon/2}$, we are done immediately because one of the
$\omega\big(n^{\ell \epsilon/2}\big)$ random $\ell$-tuples a.s. ends up in
$A_i$ and this destroys the $2$-colorability. So we can assume
that $s \leq \lfloor (1 - \epsilon/2) \log_2 n \rfloor$. Recall
that $\ell \geq 3$ and consider the following sum
$$ \sum_{s=1}^{\lfloor (1-\epsilon/2) \log_2 n \rfloor} 2^{\frac{\ell-2}{\ell-1} s}
n^{\frac{\ell}{\ell-1}(1 - \epsilon/2)}
 \leq \frac{n^{\frac{\ell-2}{\ell-1}(1 - \epsilon/2)}}{1 -
 2^{-\frac{\ell-2}{\ell-1}}} \cdot n^{\frac{\ell}{\ell-1}(1 -
 \epsilon/2)}  \leq 4 n^{2 - \epsilon}. $$
Since $G$ has at least $c n^{2-\epsilon}$ edges, there is a
subgraph $G_s$ containing at least $\frac{c}{4}
2^{\frac{\ell-2}{\ell-1} s} n^{\frac{\ell}{\ell-1}(1 -
\epsilon/2)}$ edges. As each component of $G_s$ has at most
$2^{2s}$ edges, the number of components of $G_s$ is $t =
\Omega\big(2^{-\frac{\ell}{\ell-1} s} n^{\frac{\ell}{\ell-1}(1 -
\epsilon/2)}\big)$. We set $2^s = n^{1 - \alpha}$,  $\alpha \geq
\epsilon/2$ which means that $t =
\Omega\big(n^{\frac{\ell}{\ell-1}(\alpha - \epsilon/2)}\big)$. To
summarize, we have $t$ disjoint sets $A_i$ of size
$\Theta(n^{1-\alpha})$, each of which must be monochromatic under
any feasible coloring. Thus we can apply Lemma~\ref{clusters} to
conclude that for $H= G + R$, almost surely there is no feasible
$2$-coloring.
 \hfill $\Box$

\section{General hypergraphs}

In this section we deal with the general case of a $2$-colorable
$k$-uniform hypergraph $H$, to which we add a collection of random
$\ell$-tuples $R$. Our goal is to prove Theorem~\ref{main} which
asserts that if $H$ has $\Omega\big(n^{k-\epsilon}\big)$ edges
then adding to it $\omega\big(n^{\ell \epsilon/ 2}\big)$ random
$\ell$-tuples makes it almost surely non-$2$-colorable. The proof
will proceed by induction on $k$. The base case when $k=2$ follows
from Proposition~\ref{bipartite}, so we can assume that $k
> 2$ and that the result holds for $k-1$.

We start with a series of lemmas which allow us to make
simplifying assumptions. Depending on the hypergraph $H$, we
either reduce the problem to the $(k-1)$-uniform case or prove
directly that $H + R$ is not $2$-colorable.

Since we have $\omega\big(n^{\ell \epsilon/2}\big)$ random $\ell$-tuples
available, we can divide them into a constant number of batches,
where each batch still has $\omega\big(n^{\ell \epsilon/2}\big)$
$\ell$-tuples. We will use a separate batch for each step of the
induction. We write $R = R_1 \cup R_2 \cup \ldots \cup R_k$ where
$|R_i| = \omega\big(n^{\ell \epsilon/2}\big)$ for each $i$.

\begin{lemma}
\label{k-1-degrees} Let $H_k$ be a $k$-uniform hypergraph on $n$
vertices with $c_1 n^{k-\epsilon}$ edges. Consider all
$(k-1)$-tuples $A \subset V(H_k)$ with degree greater than
$n^{1-\epsilon/2}$. If there are at least $\frac{c_1}{4}
n^{k-1-\epsilon}$ such $(k-1)$-tuples in $H_k$ then $H_k + R$ is
almost surely non-$2$-colorable.
\end{lemma}

\noindent {\bf Proof.} For each $(k-1)$-tuple $A$ of degree
$>n^{1-\epsilon/2}$, the neighborhood $N(A)$ contains
$\Omega\big(n^{\ell - \ell \epsilon/2}\big)$ distinct
$\ell$-tuples. Therefore a random $\ell$-tuple lands inside $N(A)$
with probability $\Omega\big(n^{-\ell \epsilon/2}\big)$.
Consequently, the probability that none of $\omega\big(n^{\ell
\epsilon/2}\big)$ random $\ell$-tuples from $R_k$ ends up inside
$N(A)$ is at most $\big(1-\Omega(n^{-\ell
\epsilon/2})\big)^{\omega(n^{\ell \epsilon/2})}=o(1)$. If we have
$t \geq \frac{c_1}{4} n^{k-1-\epsilon}$ such $(k-1)$-tuples, then
the expected number of them, whose neighborhood does {\em not}
contain any $\ell$-tuple in $R_k$, is $o(t)$. Therefore, by
Markov's inequality, we get almost surely at least $\frac{t}{2}
\geq \frac{c_1}{8} n^{k-1-\epsilon}$ $(k-1)$-tuples with an
$\ell$-edge in their neighborhood. Denote by $H_{k-1}$ the
$(k-1)$-uniform hypergraph formed by these $(k-1)$-tuples.

By induction, we know that $H_{k-1} + R_1 + \ldots + R_{k-1}$ is
almost surely non-$2$-colorable. Therefore for every $2$-coloring
respecting $R_1 \cup \ldots \cup R_{k-1}$, there is a
monochromatic $(k-1)$-tuple $A$ in $H_{k-1}$. Without loss of
generality assume that all vertices in $A$ are colored by $1$. By
definition, the neighborhood $N(A)$ contains an $\ell$-edge $L \in
R_k$. Either $L$ is monochromatic, or one of its vertices $x$ is
colored by $1$ as well. But then $A \cup \{x\}$ is a monochromatic
edge of $H_k$. This implies that there is no feasible $2$-coloring
for $H_k + R_1 + \ldots + R_k$. \hfill $\Box$

\vspace{0.1cm}

Thus we only need to treat the case where there are at most
$\frac{c_1}{4} n^{k-1-\epsilon}$ $(k-1)$-tuples with degree
greater than $n^{1-\epsilon/2}$, therefore at most $\frac{c_1}{4}
n^{k-\epsilon}$ edges through such $(k-1)$-tuples. We will get rid
of these high degrees by removing a constant fraction of edges and
making all degrees of $(k-1)$-tuples at most $n^{1-\epsilon/2}$.
This would also imply a bound of $n^{2-\epsilon/2}$ on the degrees
of $(k-2)$-tuples, etc. However, in the following we show that for
$(k-2)$-tuples we can assume an even stronger bound. More
specifically, we prove that if we have many edges through
$(k-2)$-tuples of degrees $n^{2-\delta}$ with $\delta \leq
\frac{\ell}{2(\ell-1)} \epsilon$, then we can proceed by
induction. For this purpose, we first show the following.

\begin{lemma}
\label{families}
Let $\ell \geq 2$ and
let $G$ be a graph on $n$ vertices with $n^{2-\delta}$ edges.
Then $G$ contains $\frac12n^{1-\delta}$ disjoint
subsets of vertices $F_1, F_2, \ldots$ such that the vertices
in each $F_j$ have disjoint neighborhoods of sizes
$d_1, d_2,\ldots$, satisfying $d_i \geq \frac12n^{1-\delta}$ and
$$ \sum{d_i^\ell} \geq \frac{n^{\ell - (\ell-1) \delta}}{2^{\ell}}\,.$$
\end{lemma}

\noindent{\bf Proof.} We iterate the following construction for
$j=1,2,\ldots,\frac12 n^{1-\delta}$.

\begin{itemize}
\item Take the vertex $v_1$ of maximum degree $d_1$ and remove
all the edges incident to its neighbors. Note that by maximality
of $d_1$, at most $d_1^2$ edges are removed.
\item In step $i$, take the vertex $v_i$ of maximum degree $d_i$
in the remaining graph and remove the edges incident to its
neighbors (again, at most $d_i^2$ edges). Repeat these steps, as
long as $\sum{d_i^2} < \frac14 n^{2-\delta}$.
\item When the procedure terminates, define $F_j = \{v_1, v_2,
 \ldots \}$. Then return to the original graph, but remove the
vertices in $F_j$ and all their edges permanently.
\end{itemize}

\begin{figure}[!ht]
\begin{center}
\scalebox{0.8}{\includegraphics{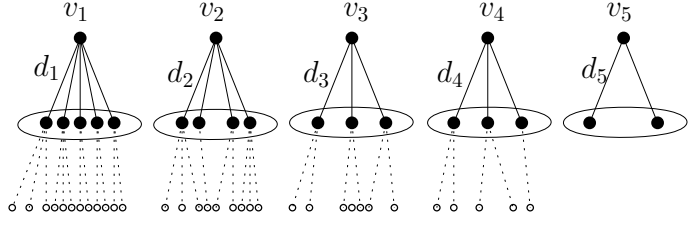}}
\caption{Construction of $F_j = \{v_1, v_2, \ldots \}$. The
neighborhood of $v_i$ is incident with at most $d_i^2$ edges.}
\end{center}
\end{figure}

By construction, the neighborhoods of the vertices in every $F_j$
are disjoint and hence with each $F_j$, we remove $\sum{d_i} \leq
n$ edges from the graph. The sets $F_j$ are also disjoint
(although the neighborhoods of vertices from different $F_j$'s are
not necessarily disjoint). Since we constructed
$\frac12n^{1-\delta}$ sets $F_j$, there are at least
$n^{2-\delta}-\frac12n^{1-\delta} \cdot n=\frac12n^{2-\delta}$
edges available at the beginning of every  construction.

Inside the construction of $F_j$, we repeat as long as $\sum{d_i^2} < \frac14
n^{2-\delta}$ and therefore we remove at most $\frac14n^{2-\delta}$
edges from the graph we started with. Hence, at every step the remaining graph still
has at least $\frac14n^{2-\delta}$ edges and so its maximum degree is at least
$\frac12n^{1-\delta}$. When we terminate we have
$\sum{d_i^2} \geq \frac14 n^{2-\delta}$. This, together with the fact that
$d_i \geq \frac12n^{1-\delta}$, implies that for every $F_j$ we have
$$\hspace{4.8cm} \sum{d_i^\ell} \geq \left( \frac12 n^{1-\delta}\right)^{\ell-2}
\sum{d_i^2} \geq \frac{n^{\ell - (\ell-1) \delta}}{2^{\ell}}.
\hspace{4.8cm} \Box$$

\begin{lemma}
\label{k-2-degrees} Let $H_k$ be a $k$-uniform hypergraph on $n$
vertices with $c_1 n^{k-\epsilon}$ edges. Consider $(k-2)$-tuples
of degree $n^{2-\delta}$ where $\delta \leq \frac{\ell
}{2(\ell-1)} \epsilon$. If there are at least $\frac{c_1}{4}
n^{k-\epsilon}$ edges through such $(k-2)$-tuples then $H_k + R$
is almost surely non-$2$-colorable.
\end{lemma}

\noindent {\bf Proof.} Consider a $(k-2)$-tuple $A$ of degree
$n^{2-\delta}$. The link of $A$ in $H_k$ is a graph $\Gamma(A)$
with $n^{2-\delta}$ edges. By Lemma~\ref{families}, we find
$\frac12 n^{1-\delta}$ subsets $F_j$ such that vertices in $F_j$
have disjoint neighborhoods in $\Gamma(A)$ with sizes satisfying
$\sum{d_i^\ell} \geq 2^{-\ell} n^{\ell - (\ell-1) \delta}$. We
repeat this construction for each $(k-2)$-tuple of degree
$n^{2-\delta}$ with $\delta \leq \frac{\ell}{2(\ell-1)} \epsilon$.
For each of them, we construct $\frac12 n^{1-\delta}$ sets as
above. Assuming that the total number of edges through such
$(k-2)$-tuples is at least $\frac{c_1}{4} n^{k-\epsilon}$, we get
$\frac{c_1}{8} n^{k-1-\epsilon}$ sets $F_j$ in total.

Now fix a set $F_j$. Call it {\em good} if after adding random
$\ell$-tuples from $R_k$ there is at least one vertex in $F_j$
whose neighborhood in $\Gamma(A)$ contains a random $\ell$-tuple.
If this is not the case, call it {\em bad}. We estimate the
probability that $F_j$ is bad. By Lemma~\ref{families}, the total
number of $\ell$-tuples in the neighborhoods of vertices in $F_j$
is
$$ \sum{d_i \choose \ell}=\Omega \left(\sum\frac{d_i^\ell}{\ell!}\right) =
\Omega \left(
 \frac{n^{\ell - (\ell-1) \delta}}{2^\ell \ell!}\right)=
\Omega \big( n^{\ell - \ell \epsilon/2}\big). $$ Thus the
probability that a random $\ell$-tuple falls inside some
neighborhood of $F_j$ is $\sum{d_i \choose \ell} / {n \choose
\ell} = \Omega \big( n^{-\ell \epsilon / 2}\big)$. After adding
the entire batch of random $\ell$-tuples $R_k$,
$$ \Pr\big[F_j \mbox{ is bad}\big]
 = \left(1 - \Omega\big( n^{-\ell \epsilon / 2}\big)\right)^{-\omega(n^{\ell
 \epsilon / 2})}= o(1).$$
Consequently, the expected fraction of bad $F_j$'s is $o(1)$. By
Markov's inequality, this fraction is almost surely at most one
half, which means that at least $\frac{c_1}{16} n^{k-1-\epsilon}$
sets $F_j$ have a vertex $v \in F_j$ whose neighborhood contains
some $\ell$-tuple from $R_k$. For each such $F_j$, we have a set
$A$ of size $k-2$ which together with $v$ forms a $(k-1)$-tuple
whose neighborhood in $H_k$ contains an $\ell$-tuple from $R_k$.
We could get the same $(k-1)$-tuple in $k-1$ different ways, but
in any case we have at least $\frac{c_1}{16k} n^{k-1-\epsilon}$
such $(k-1)$-tuples which form an edge set of a $(k-1)$-uniform
hypergraph $H_{k-1}$.

By the induction hypothesis, $H_{k-1} + R_1 + \ldots + R_{k-1}$ is
almost surely non-$2$-colorable. Therefore, for any $2$-coloring
which respects the $\ell$-edges from $R_1 + \ldots + R_{k-1}$,
there must be a monochromatic $(k-1)$-edge $B$ in $H_{k-1}$.
However, since there is an $\ell$-edge from $R_k$ in the
neighborhood of $B$, one of its vertices should have the same
color as $B$. This would form a monochromatic edge in $H_k$ so
there is no feasible $2$-coloring for $H_k + R_1 + \ldots + R_k$.
\hfill $ \Box$

Thus we can also assume that at most $\frac{c_1}{4}
n^{k-\epsilon}$ edges go through $(k-2)$-tuples of degree
$n^{2-\delta}, \delta \leq \frac{\ell}{2(\ell-1)} \epsilon$.
Before the last part of the proof, we make further restrictions on
the degree bounds and structure of our hypergraph, by finding a
subhypergraph $H_\alpha$ described in the following lemma.

\begin{lemma}
\label{k-partite} Let $H_k = (V,E)$ be a $k$-uniform hypergraph
with $c_1 n^{k-\epsilon}$ edges, such that at most $\frac{c_1}{4}
n^{k-\epsilon}$ edges go through $(k-1)$-tuples of degree $\geq
n^{1-\epsilon/2}$ and at most $\frac{c_1}{4} n^{k-\epsilon}$ edges
go through $(k-2)$-tuples of degree $n^{2-\delta}, \delta \leq
\frac{\ell}{2(\ell-1)} \epsilon$. Then for some constant $\alpha
\geq \epsilon/2$, $H_k$ contains a subhypergraph $H_\alpha$ with
the following properties
\begin{enumerate}
 \item $H_\alpha$ is $k$-partite, i.e. $V$ can be partitioned
 into $V_1 \cup V_2 \cup \ldots \cup V_k$ so that every edge of
 $H_\alpha$ intersects each $V_i$ in one vertex.
 \item Every vertex has degree at most $n^{k-1 - \frac{\ell}
 {2(\ell-1)}  \epsilon}$.
 \item The degree of every $(k-1)$-tuple in $V_1 \times V_2 \times \ldots \times V_{k-1}$
is either $0$ or between $n^{1-\alpha}$ and $2n^{1-\alpha}$.
 \item The number of edges in $H_\alpha$ is at least
 $$ c_5 \left(n^{k - \epsilon - \frac{\epsilon-\alpha}{\ell-1}} +
 n^{k - \epsilon - \frac{\ell-2}{\ell-1} (\alpha - \epsilon/2)}\right), $$
for some constant $c_5 =c_5(k,\ell,c_1)$.
\end{enumerate}
\end{lemma}

\noindent {\bf Proof.} First, remove all edges through
$(k-1)$-tuples of degree $\geq n^{1-\epsilon/2}$ and through
$(k-2)$-tuples of degree $n^{2-\delta}, \delta \leq
\frac{\ell}{2(\ell-1)} \epsilon$. We get a hypergraph $H'$ such
that the degrees of all $(k-1)$-tuples are at most
$n^{1-\epsilon/2}$, the degrees of all $(k-2)$-tuples are at most
$n^{2 - \frac{\ell}{2(\ell-1)} \epsilon}$, and the number of edges
is at least $c_2 n^{k-\epsilon}$ edges, $c_2 = c_1/2$.
Consequently, the degree of every vertex in $H'$ is at most
$n^{k-3} \cdot n^{2 - \frac{\ell}{2(\ell-1)}  \epsilon} = n^{k-1 -
\frac{\ell}{2(\ell-1)}  \epsilon}$.

Next, we use a well known fact, proved by Erd\H{o}s and Kleitman
\cite{EK} that every $k$-uniform hypergraph $H'$ with $c_2
n^{k-\epsilon}$ edges contains a $k$-partite subhypergraph with at
least $c_3 n^{k-\epsilon}$ edges where $c_3 = \frac{k!}{k^k}c_2$.
This can be achieved for example by taking a random partition of
the vertex set into $k$ parts and computing the expected number of
edges which intersect all of them. Let $(V_1, V_2, \ldots, V_k)$
be a partition, so that at least $c_3 n^{k-\epsilon}$ edges of
$H'$ have one vertex in every $V_i$. Discard all other edges and
denote this $k$-partite hypergraph by $H''$.

Consider all $(k-1)$-tuples in $V_1 \times V_2 \times \ldots
\times V_{k-1}$ whose degree in $H''$ is less than $\frac{c_3}{2}
n^{1-\epsilon}$. Delete all their edges, which is at most ${n
\choose k-1} \frac{c_3}{2} n^{1-\epsilon} \leq \frac{c_3}{2}
n^{k-\epsilon}$ edges in total. We still have at least $c_4
n^{k-\epsilon}$ edges, where $c_4 = c_3/2$. Now the degree of
every $(k-1)$-tuple in $V_1 \times V_2 \times \ldots \times
V_{k-1}$ is either $0$ or between $c_4 n^{1-\epsilon}$ and
$n^{1-\epsilon/2}$. Finally, we are going to find a subhypergraph
in which all the non-zero degrees of $(k-1)$-tuples are
$\Theta(n^{1-\alpha})$ and the number of edges is at least
$$ c_5 \left(n^{k - \epsilon - \frac{\epsilon-\alpha}{\ell-1}} +
 n^{k - \epsilon - \frac{\ell-2}{\ell-1} (\alpha - \epsilon/2)}\right). $$
The existence of such a subhypergraph can be proved by an
elementary counting argument. Let $n^{1-\alpha} = 2^i$ and
partition $V_1 \times V_2 \times \ldots \times V_{k-1}$ into
groups of $(k-1)$-tuples with degrees in intervals $[2^i,
2^{i+1})$, where $i$ ranging between $i_1 = \log_2 (c_4
n^{1-\epsilon})$ and $i_2 = \log_2 (n^{1 - \epsilon/2})$. Consider
the following two expressions:
$$\sum_{i=i_1}^{i_2} {2^{-i/(\ell-1)}}
\leq \frac{(c_4
n^{1-\epsilon})^{-\frac{1}{\ell-1}}}{1-2^{-\frac{1}{\ell-1}}}
 \leq  2(\ell-1) c_4^{-1} n^{-\frac{1-\epsilon}{\ell-1}}$$
and
$$\sum_{i=i_1}^{i_2} {2^{\frac{\ell-2}{\ell-1} i}} \leq
\frac{n^{\frac{\ell-2}{\ell-1} (1 - \epsilon/2)}}{1 -
2^{-\frac{\ell-2}{\ell-1}}}
 \leq 4 n^{\frac{\ell-2}{\ell-1} (1 - \epsilon/2)}.$$
Normalizing by the right-hand side and taking the average, we get
$$ \sum_{i=i_1}^{i_2} \left( \frac{2^{-\frac{i}{\ell-1}}}
{4(\ell-1)c_4^{-1} n^{-\frac{1-\epsilon}{\ell-1}}} +
\frac{2^{\frac{\ell-2}{\ell-1}i}}{8
n^{\frac{\ell-2}{\ell-1}(1-\epsilon/2)}} \right) \leq 1$$ By the
pigeonhole principle, there is an $i$ such that the fraction of
edges through $(k-1)$-tuples with degree between
$2^i=n^{1-\alpha}$ and $2^{i+1}=2n^{1-\alpha}$ is at least
$$ \frac{2^{-\frac{i}{\ell-1}}}{4(\ell-1)c_4^{-1} n^{-\frac{1-\epsilon}{\ell-1}}} +
\frac{2^{\frac{\ell-2}{\ell-1}i}}{8 n^{\frac{\ell-2}{\ell-1}
(1-\epsilon/2)}} = \frac{c_4}{4 (\ell-1)}
n^{-\frac{\epsilon-\alpha}{\ell-1}} + \frac{1}{8}
n^{-\frac{\ell-2}{\ell-1} (\alpha - \epsilon/2)} $$ so the lemma
holds with $c_5 = c_4 \cdot \min\big\{\frac{c_4}{4 (\ell-1)},
\frac{1}{8} \big\}$. \hfill $\Box$

\vspace{0.1cm}

Note that in this lemma, we lose more than a constant fraction of
the edges. However, from now on, we do not use induction anymore
and will prove directly that $H_\alpha + R$ is almost surely
non-$2$-colorable. We will proceed in $t = c_5 \ell^{-k}
n^{\frac{\ell}{\ell-1} (\alpha - \epsilon/2)}$ stages. For each
stage, we allocate a certain number of random $\ell$-tuples.
Namely, we set again $R = R_1 \cup R_2 \cup \ldots \cup R_k,\,
|R_i|=\omega\big(n^{\ell \epsilon / 2}\big)$. Furthermore, we
divide each $R_j$ for $j \leq k-1$ into $t$ parts $R_{1,j},
\ldots, R_{t,j}$ so that
 $$ |R_{i,j}| = \omega\left(\frac{n^{\ell \epsilon / 2}}{t}\right)
  = \omega\left(n^{\ell \epsilon/2 - \frac{\ell}{\ell-1} (\alpha -
  \epsilon/2)}\right).$$
The random set $R_{i,j}$ will be used for the $j$-th ``level" of
the $i$-th stage. The following lemma describes one stage of the
construction. Finally, $R_k$ will be used in the last step of the proof.

\begin{lemma}
\label{one-tree} Let $H_\alpha$ be a $k$-uniform $k$-partite
hypergraph where the degree of every $(k-1)$-tuple in $V_1 \times
V_2 \times \ldots \times V_{k-1}$ is either zero or
is in the interval $[n^{1-\alpha}, 2n^{1-\alpha}]$,
and the number of edges in $H_\alpha$ is at least
$$ c_5 n^{k - \epsilon - \frac{\ell-2}{\ell-1} (\alpha - \epsilon/2)}.$$
Then almost surely, there exists a family of $q = \ell^{k-2}$ sets
$S_1, \ldots, S_q$, $n^{1-\alpha} \leq S_i \leq 2n^{1-\alpha}$,
such that for every feasible $2$-coloring of $H_\alpha + R_{i,1} +
\ldots + R_{i,k-1}$ at least one $S_i$ is monochromatic.
\end{lemma}

\noindent {\bf Proof.} We are going to construct an $\ell$-ary
tree $T$ of depth $k-1$. We denote vertices on the $j$-th level by
$v_{a_1 a_2 \ldots a_{j-1}}$ where $a_i \in \{1,2,\ldots,\ell\}$.
$T$ is rooted at a vertex in $V_1$ and the $j$-th level is
contained in $V_j$. We construct $T$ in such a way that the
vertices along every path which starts at the root and has length
$k-1$ form a $(k-1)$-tuple with degree $\Theta(n^{1-\alpha})$ in
$H_\alpha$. The neighborhoods of all branches of length $k-1$ will
be our sets $S_i$ (not necessarily disjoint). In addition, the set
of $\ell$ children of every node on each level $j \leq k-2$, like
$\{v_{a_1 a_2 \ldots a_{j-1} 1}, v_{a_1 a_2 \ldots a_{j-1}
2},\ldots, v_{a_1 a_2 \ldots a_{j-1} \ell}\}$, will form an edge
of $R_{i,j}$.

\begin{figure}[!ht]
\begin{center}
\scalebox{0.8}{\includegraphics{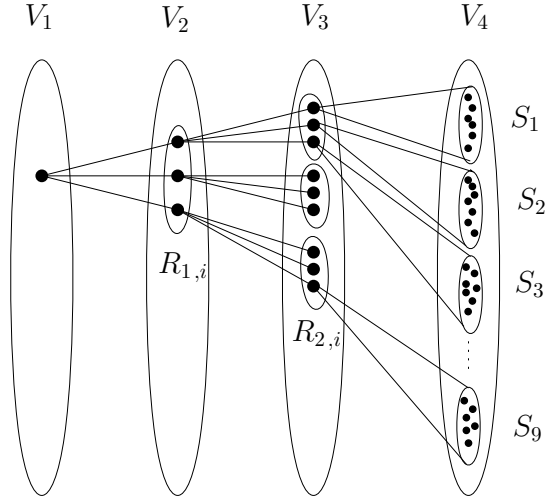}} \caption{Construction
of the tree $T$, for $k=4$ and $\ell=3$. Branches of the tree form
active $(k-1)$-tuples, with neighborhoods $S_i$. Each set of
children on level $j+1$ forms an edge of $R_{i,j}$.}
\end{center}
\end{figure}

Assuming the existence of such a tree, consider any $2$-coloring
of $H_\alpha + R_{i,1} + \ldots + R_{i,k-1}$. Since the children
of each vertex on level $j < k-1$ form an $\ell$-edge in
$R_{i,j}$, every vertex has children of both colors. In
particular, there is always one child with the same color as its
parent. Therefore, starting from the root, we can always find a
monochromatic branch $A$ of length $k-1$. Since all the extensions
of this branch to edges of $H_\alpha$ must be $2$-colored, all the
vertices in $S_i = N(A)$ must have the same color.

We grow the tree level by level, maintaining the property that all
branches have sufficiently many extensions to edges of $H_\alpha$.
More precisely, we call an $r$-tuple in $V_1 \times \ldots \times
V_r$ {\em active} if its degree is at least
 $$ \Delta_r = \frac{c_5}{2^r} n^{k-r-\epsilon - \frac{\ell-2}{\ell-1}(\alpha -
 \epsilon/2)}.$$

\noindent {\bf Claim.} Every active $r$-tuple $A$, $r \leq k-2$,
can be extended to at least
 $$ d_r = \frac{\Delta_r}{4 n^{k-r-1-\alpha}} = \frac{c_5}{2^{r+2}}
 n^{1 - \epsilon/2 + \frac{1}{\ell-1} (\alpha - \epsilon/2)} $$ active
$(r+1)$-tuples $A \cup \{x\}, x \in V_{r+1}$.

{\em Proof.} Suppose that fewer than $d_r$ extensions of $A$ are
active. Since the degrees of $(k-1)$-tuples are at most $2
n^{1-\alpha}$, we get that any $(r+1)$-tuple has degree at most $2
n^{k-r-1-\alpha}$. Therefore the number of edges through all
active extensions of $A$ is smaller than $d_r \cdot 2
n^{k-r-1-\alpha}=\frac12\Delta_r$. We also have inactive
extensions of $A$ which have degrees less than $\Delta_{r+1}$. The
total number of edges through these extensions of $A$ is smaller
than $n \Delta_{r+1} = \frac12 \Delta_r$. But the total number of
edges through $A$ is at least $\Delta_r$. This contradiction
proves the claim. \hfill $\Box$

We start our construction from an active vertex $v \in V_1$. Since
$H_\alpha$ has at least $n \Delta_1$ edges, such a vertex must
exist. By our claim, $v$ can be extended to at least $d_1$ active
pairs $\{v,x\}, x \in W_2 \subset V_2$. Consider this set of $d_1$
vertices $W_2$. The probability that a random $\ell$-tuple falls
inside $W_2$ is ${d_1 \choose \ell} / {n \choose \ell} =
\Omega(n^{-\ell \epsilon/2 + \frac{\ell}{\ell-1} (\alpha -
\epsilon/2)})$. Now we use $\omega(n^{\ell \epsilon/2  -
\frac{\ell}{\ell-1} (\alpha - \epsilon/2)})$ random $\ell$-tuples
from $R_{i,1}$ that we allocated for the first level of this
construction. This means that almost surely, we get an $\ell$-edge
$\{v_1, \ldots, v_{\ell}\} \in R_{i,1}$ such that $\{v, v_i\}$ is
an active pair for each $i=1,2,\ldots,\ell$.

We continue growing the tree, using the random $\ell$-tuples of
$R_{i,j}$ on level $j$. Since we have ensured that each path from
the root to the level $j$ from an active $j$-tuple, it has at
least $d_j$ extensions to an active $(j+1)$-tuple. Again, the
probability that a random $\ell$-tuple hits the extension vertices
$W_{j+1} \subset V_{j+1}$ for a given path is ${d_j \choose \ell}
/ {n \choose \ell} = \Omega\big(n^{-\ell \epsilon/2 +
\frac{\ell}{\ell-1} (\alpha - \epsilon/2)}\big)$. Almost surely,
one of the $\ell$-tuples in $R_{i,j}$ will hit these extension
vertices and we can extend this path to $\ell$ children on level
$j+1$. The number of paths from the root to level $j$ is bounded
by $\ell^{j-1}$ which is a constant, so in fact we will almost
surely succeed to build the entire level.

In this way, we a.s. build the tree all the way to level $k-1$.
Every path from the root to one of the leaves forms an active
$(k-1)$-tuple and has degree $\in [n^{1-\alpha}, 2n^{1-\alpha}]$.
Define $S_1, S_2, \ldots, S_q$ to be the neighborhoods of all
these $q = \ell^{k-2}$ paths. By construction, for any feasible
$2$-coloring of $H_\alpha + R_{i,1} + \ldots + R_{i,k-1}$, one of
these paths is monochromatic which implies that the corresponding
set $S_i$ is monochromatic as well. \hfill $\Box$

\begin{lemma}
\label{trees} Let $H_\alpha$ be a $k$-uniform $k$-partite
hypergraph where the degree of every vertex is at most
$n^{k-1-\frac{\ell}{2(\ell-1)} \epsilon}$, the degree of every
$(k-1)$-tuple in $V_1 \times V_2 \times \ldots \times V_{k-1}$ is
either zero or is in the interval $[n^{1-\alpha}, 2n^{1-\alpha}]$, and the number of
edges in $H_\alpha$ is at least
$$ c_5 n^{k - \epsilon - \frac{\epsilon-\alpha}{\ell-1}} +
 c_5 n^{k - \epsilon - \frac{\ell-2}{\ell-1} (\alpha - \epsilon/2)}.$$
Then almost surely, $H_\alpha + R$ is not $2$-colorable.
\end{lemma}

\noindent {\bf Proof.} We apply Lemma~\ref{one-tree} repeatedly in
$t = c_5 \ell^{-k} n^{\frac{\ell}{\ell-1} (\alpha - \epsilon/2)}$
stages. In each stage $i$, we almost surely obtain $q=\ell^{k-2}$
sets $S_{i,1}, \ldots, S_{i,q}$, $n^{1-\alpha} \leq |S_{i,j}| \leq
2n^{1-\alpha}$ such that for any $2$-coloring of the hypergraph
$H_\alpha+\sum R_{i,j}$, one of these sets must be monochromatic.
If this happens, we call such a stage ``successful". After each
successful stage, we remove all edges of $H_\alpha$ incident with
any of the sets $S_{i,1}, \ldots, S_{i,q}$. Since degrees are
bounded by $n^{k-1-\frac{\ell}{2(\ell-1)} \epsilon}$ and we repeat
$t = c_5 \ell^{-k}  n^{\frac{\ell}{\ell-1} (\alpha - \epsilon/2)}$
times, the total number of edges we remove is at most
$$ \sum_{i=1}^{t} \sum_{j=1}^{q} |S_{i,j}| n^{k-1-\frac{\ell}{2(\ell-1)} \epsilon}
 \leq t q \cdot 2 n^{1-\alpha} \cdot n^{k-1-\frac{\ell}{2(\ell-1)} \epsilon} =
2 c_5 \ell^{-2}n^{k-\epsilon - \frac{\epsilon - \alpha}{\ell-1}}
\leq c_5 n^{k - \epsilon - \frac{\epsilon - \alpha}{\ell-1}}. $$
In particular, before every stage we still have at least $c_5 n^{k
- \epsilon - \frac{\ell-2}{\ell-1} (\alpha - \epsilon/2)}$ edges
available, so we can use Lemma~\ref{one-tree}. Since the expected
number of stages that are not successful is $o(t)$, by Markov's
inequality, we almost surely get at least $t/2$ successful stages.
Eventually, we obtain sets $S_{i,j}$ for $1 \leq i \leq t/2$ and
$1 \leq j \leq q$ such that
\begin{itemize}
\item For $i_1 \neq i_2$ and any $j_1, j_2$, $S_{i_1, j_1} \cap
S_{i_2, j_2} = \emptyset$.
\item For any $2$-coloring of $H_\alpha+\sum R_{i,j}$ and any $i$, there is $j_i$ such
that
$S_{i,j_i}$ is monochromatic.
\end{itemize}

Finally, we add once again a collection $R_k$ of $\omega(n^{\ell
\epsilon/2})$ random $\ell$-tuples. We do not know a priori
which selection of sets $S_{i,j}$ will be monochromatic but there is only
exponential number of choices $q^{t/2} =
e^{O(t)}$. For any specific choice of sets to be monochromatic,
Lemma~\ref{clusters} says that the probability that after adding
$\omega(n^{\ell \epsilon/2})$ random $\ell$-tuples, there is a feasible
$2$-coloring keeping these sets monochromatic, is $e^{-\omega(t)}$.
By the union bound, the probability that there exist a proper $2$-coloring
of $H_\alpha+\sum R_{i,j}+R_k$ is at most $q^{t/2}e^{-\omega(t)}=o(1)$.
This completes the proof of this lemma together with the proof of Theorem~\ref{main}.
\hfill $\Box$

\vspace{0.4cm}
\noindent
{\bf  Acknowledgment.}\, The first author would like to thank Michael
Krivelevich for helpful and stimulating discussions.

\end{document}